\documentclass[12pt,oneside]{article}
\usepackage{latexsym}
\usepackage{amsfonts}
\usepackage{amsmath, amsthm, amssymb}
\usepackage{setspace}
\usepackage{graphicx}

\newcommand{\SX}{S(X,\varphi,Q)}

\newcommand{\UX}{U(X,\varphi,P)}

\newcommand{\GsX}{G^s(X,\varphi,Q)}
\newcommand{\GuX}{G^u(X,\varphi,P)}

\newcommand{\GsSig}{G^s(\Sigma,\sigma,Q)}
\newcommand{\GuSig}{G^u(\Sigma,\sigma,P)}

\newtheorem{theorem}{Theorem}[section]
\newtheorem{proposition}[theorem]{Proposition}

\newtheorem{lemma}[theorem]{Lemma}
\newtheorem{definition}[theorem]{Definition}

\title{Trace Asymptotics for $C^{\ast}$-Algebras from Smale Spaces}
\author{D.B. KILLOUGH, \\
Department of Mathematics, Physics, and Engineering,\\
Mount Royal University,\\
Calgary, AB, Canada T3E 6K6,\\
bkillough@mtroyal.ca\\ 
I.F. PUTNAM\footnote{Supported in part by an NSERC Discovery Grant}, \\
Department of Mathematics and Statistics,\\
University of Victoria,\\
Victoria, B.C., Canada V8W 3R4,\\
ifputnam@uvic.ca\\
}
\date{\today}
\begin{document}

\maketitle

\abstract
We consider $C^{\ast}$-algebras associated with stable and unstable equivalence in hyperbolic dynamical systems known as Smale spaces. These systems
include shifts of finite type, in which case these $C^{*}$-algebras are both AF-algebras. These algebras have  fundamental
representations on a single Hilbert space (subject to a choice of periodic points) which have a number of special properties. In particular,
the product between any element of the first algebra with one from the second is compact. In addition, there is a single
unitary operator which implements actions on both. Here, under the hypothesis that the system is mixing, we show that the 
(semi-finite) traces on these algebras may be obtained through a limiting process and the usual operator trace.

\section{Introduction}

\subsection{Smale Spaces}
A Smale space, as defined by David Ruelle \cite{ruelle78}, is a compact metric space, $X$, together
with a homeomorphism, $\varphi$, which is hyperbolic. 
These include the basic sets of
Smale's Axiom A systems \cite{smale67}. Another special case of great interest are the shifts of finite type \cite{bowen78, lindmarcus}.

Informally, the structure of $(X,\varphi)$ is such that, for each point $x$ in $X$ and each $\epsilon >0 $ 
and sufficiently small, there are sets $X^s(x,\epsilon)$ and $X^u(x,\epsilon)$, called the local stable and unstable sets respectively. Their Cartesian product is homeomorphic to a neighbourhood of $x$ (in a canonical and 
$\varphi$-invariant way)
and the map $\varphi$ is (uniformly) contracting on the former while $\varphi^{-1}$ is contracting on the latter.

The precise  axiom for a Smale space is the existence of a map defined on pairs $(x,y)$ in $X \times X$ which are sufficiently close. We have a constant $\epsilon_{X} > 0$ and for $x, y$ in $X$ with 
$d(x, y) \leq \epsilon_{X}$, the image of $(x, y)$ is denoted
$[x, y]$.  
This map satisfies a number of identities and we refer the reader to \cite{ruelle78} or \cite{putnam96} for a complete discussion.
From the map $[,]$, we then define $X^{s}(x, \epsilon) = \{ y \in X \mid d(x,y) < \epsilon, [y, x] = x \}$
and $X^{u}(x, \epsilon) = \{ y \in X \mid d(x,y) < \epsilon, [x,y] = x \}$, for any $0 < \epsilon \leq \epsilon_{X}$.
The contracting and expanding properties of the map $\varphi$ on $X^s(x,\epsilon)$ and $X^u(x,\epsilon)$
are that there is a constant $0 < k_e < 1$ such that
\begin{eqnarray*}
 d(\varphi(y), \varphi(z)) & \leq  & k_e d(y, z), y, z \in X^s(x,\epsilon), \\
  d(\varphi^{-1}(y), \varphi^{-1}(z)) & \leq  & k_e d(y, z), y, z \in X^u(x,\epsilon),
\end{eqnarray*}
which are axioms.

There is also a notion of a global stable (unstable) set for a point $x$, which we denote $X^s(x)$ ($X^u(x)$).  
This is simply the set of all points $y \in X$ such that $d(\varphi^{n}(x), \varphi^{n}(y)) \rightarrow 0$ as $n \rightarrow +\infty \ (-\infty)$.  
The collection of sets $\{ X^s(y,\delta) \ | \ y \in X^s(x), \ \delta > 0 \}$ forms a neighbourhood base for a topology on 
$X^s(x)$ in which it is locally compact and Hausdorff.  This is the topology that we use on $X^s(x)$ (not the relative topology from $X$). 
 There is an analogous topology on $X^u(x)$. Briefly, we have two equivalence relations on $X$ which are
 transverse in some sense. It is this transversality which lies at the core of the results in \cite{kpw} and
 our main result below.
 
 We say that $(X, \varphi)$ is irreducible if, for every (ordered) pair of non-empty open sets $U$ and $V$,
 there exists an $N \geq 1$ such that $\varphi^{N}(U) \cap V \neq \emptyset$. Also, we
 say that $(X, \varphi)$ is mixing if, for every (ordered) pair of non-empty open sets $U$ and $V$,
 there exists an $N_{0} \geq 1$ such that $\varphi^{N}(U) \cap V \neq \emptyset$, for all $N \geq N_{0}$.

An important feature of an irreducible Smale space $(X,\varphi)$ is the existence and uniqueness 
of a  $\varphi$-invariant probability measure maximizing the entropy 
of $\varphi$, see \cite{ruellesullivan, katokhass}.  We call this the Bowen 
measure and denote it by $\mu_X$, or when the space is obvious, simply $\mu$.  
%In \cite{}, Bowen constructed this measure as a limit of measures supported on periodic points, and in \cite{killoughputnam11} it was shown that a similar construction could be carried out with heteroclinic points.  In the case that the Smale space is a SFT, the Bowen measure is usually referred to as the Parry measure.
It was shown in \cite{ruellesullivan} (but also see \cite{killoughputnam11} for further discussion) that, 
as our space is locally a product space, the 
Bowen measure is locally a product measure. 
 Specifically, for each $x$ in $X$, we  have measures $\mu^{s,x}$ and $\mu^{u,x}$ 
defined on $X^{s}(x)$ and $X^{u}(x)$, respectively.
Secondly, the measure $\mu^{s,x}$  depends only on the stable 
equivalence class of $x$; that is, if $y$ is in $X^{s}(x)$, then 
$\mu^{s,x} = \mu^{s,y}$. (Put another way, we should be writing 
$\mu^{s,X^{s}(x)}$, but that notation is rather too clumsy.)
A similar statement holds for $\mu^{u,x}$. 
It is worth noting that these measures are not finite, but are regular Borel measures.
Moreover, these satisfy the following conditions.

\begin{theorem}
\begin{enumerate}
\item For all $x$ in $X$, $\epsilon > 0$ and Borel sets 
$B \subset X^{u}(x, \epsilon)$ and $C \subset X^{s}(x, \epsilon)$, we have 
\[
 \mu([B, C]) = \mu^{u,x}(B)\mu^{s,x}(C)
\]
whenever $\epsilon$ is sufficiently small so that $[B, C]$ is defined.
\item For $x$, $y$ in $X$, $\epsilon > 0$ and  a Borel set $B \subset X^{u}(x, \epsilon)$, we have 
\[
 \mu^{u,y}([B, y])  = \mu^{u,x}(B),
\]
whenever $d(x,y)$ and $\epsilon$ are sufficiently small so that $[B, y]$ is defined.
\item For $x$, $y$ in $X$, $\epsilon > 0$ and  a Borel set $C \subset X^{s}(x, \epsilon)$, we have 
\[
 \mu^{s,y}([y, C])  = \mu^{s,x}(C),
\]
whenever $d(x,y)$ and $\epsilon$ are sufficiently small so that $[y, C]$ is defined.
\item $\mu^{s,\varphi(x)} \circ \varphi = \lambda^{-1} \mu^{s,x}$.
\item $\mu^{u,\varphi(x)} \circ \varphi = \lambda \mu^{u,x}$.
\end{enumerate}
Here $\log(\lambda)$ is the topological entropy of $(X,\varphi)$.
\end{theorem}

In future, it will cause no confusion to drop the superscript $x$.
In the case that the Smale space is  a shift of finite type (SFT), the Bowen measure is the same as the Parry measure. 
%We  present a brief description of the Parry measure for a \emph{mixing} SFT in section 3.

\subsection{$C^{*}$-algebras}

Several $C^{\ast}$-algebras can be constructed from a given Smale space. 
The construction is originally due to Ruelle \cite{ruelle88}. In the case of a 
shift of finite type, these  algebras were defined and studied 
earlier by Krieger, Cuntz and Krieger in \cite{cuntzkrieger,  krieger80}.  
In the more general Smale space setting we refer the reader to \cite{putnam96, putnam00, putnam05}.  The algebras that we will be concerned with for a Smale space $(X, \varphi)$ are known as the stable and unstable algebras.   

We consider the groupoid $C^{*}$-algebras associated with stable and unstable equivalence.
Here, it is most convenient to consider stable equivalence on unstable sets, since these
function as abstract transversals. Also, it is convenient to have such sets which are $\varphi$-invariant.
So we fix a finite $\varphi$-invariant set, $P$, (the periodic points in an irreducible Smale 
space are always dense, so there is a good supply of such sets) and let
$$
X^s(P) = \bigcup_{p \in P} X^s(p), \ X^u(P) = \bigcup_{p \in P} X^u(p),
$$
We then define
\begin{eqnarray*}
%G^h(X) &=& \{(x,y) \ | \ x \stackrel{h}{\sim} y \} \\
G^s(X,\varphi,P) &=& \{(x,y) \ | \ x \stackrel{s}{\sim} y,  \ x, y \in X^u(P) \} \\
G^u(X,\varphi,P) &=& \{(x,y) \ | \ x \stackrel{u}{\sim} y,  \ x, y \in X^s(P) \}
\end{eqnarray*}
These groupoids have topologies, which will be described in detail in section 2,
in which they are \'{e}tale. We define $S(X, \varphi, P) = C^{*}(G^s(X, \varphi, P))$ and 
$U(X, \varphi, P) = C^{*}(G^u(X, \varphi, P))$.

We summarize some of the properties of $S(X, \varphi, P)$ and 
$U(X, \varphi, P)$. Both are simple if and only if $(X, \varphi)$ is mixing. They are both amenable, stable
and finite. In the case of a shift of finite type, each is an AF-algebra. We also point out that the homeomorphism $\varphi$ yields an automorphism, $\alpha$, of each of these $C^*$-algebras: $\alpha(a)(x,y) = a(\varphi^{-1}(x),\varphi^{-1}(y))$.

Our system of measures $\mu^{u,x}$ and $\mu^{s,x}$ provide faithful, semi-finite traces on 
$S(X, \varphi, P)$ and 
$U(X, \varphi, P)$ as follows: $\sum_{x \in P} \mu^{u, x}$ defines a measure on 
the unit space of $G^s(X,\varphi,P)$ and the invariance property 2 of Theorem
1.1 implies that this will be a trace.
Again, we will denote this measure simply by $\mu^{u}$.
We summarize the following result. A version first appears in \cite{putnam96}. 
Our definition of $\SX$, which differs from that in \cite{putnam96}, can be found in \cite{putnam00, putnamspielberg, killoughputnam10}, however, the proofs used in \cite{putnam96} apply with little alteration.

\begin{theorem}
Let $(X, \varphi)$ be an irreducible Smale space and let $P,Q$
be finite $\varphi$-invariant subsets of $X$.
For $a \in C_c(G^{s}(X, \varphi,Q))$ and $b \in C_c(G^{u}(X, \varphi, P))$, define
\begin{eqnarray*}
% \tau^h(f) &=& \int_X f(x,x) d\mu, \\
\tau^s(a) &=&  \int_{X^u(Q)} a(x,x) d\mu^u, \ \textrm{and} \\
\tau^u(b) &=& \int_{X^s(P)} b(x,x) d\mu^s.
\end{eqnarray*}
%$\tau^h$ extends to a bounded trace on $\HX$, while 
Then $\tau^s$ and $\tau^u$ extend to semi-finite traces on $\SX$ and $\UX$. 
\end{theorem}

\subsection{Fundamental representation}

Here, we choose two finite $\varphi$-invariant sets, $P$ and $Q$. We will consider the
$C^{*}$-algebras $\SX$ and $\UX$. Define
\[
 X^{h}(P, Q) = X^{s}(P) \cap X^{u}(Q).
\]
(The `h' stands for heteroclinic.)  This is a subset of $X^{u}(Q)$, the unit space of  $G^{s}(X, \varphi, Q)$
and consists of some stable equivalence classes in that set. Hence, the Hilbert space
$l^{2}(X^{h}(P,Q))$ admits a natural representation of $S(X, \varphi, Q)$ which 
we will describe in detail in the next section. In an analogous way, it also
admits a representation of $U(X, \varphi, P)$.
We suppress these representations in our notation and
simply assume that these $C^{*}$-algebras are acting on this Hilbert space.  We also note that $X^{h}(P, Q)$ is countable
and that $\delta_{x}, x \in X^{h}(P, Q)$ denotes the usual orthonormal basis.
We note that $u \xi = \xi \circ \varphi^{-1}$ is a unitary operator on the 
Hilbert space and implements the automorphisms $\alpha$ of $S(X, \varphi, Q)$ and $U(X, \varphi, P)$.

This pair of $C^{*}$-algebras acting on this Hilbert space possesses a number of interesting features. 
We summarize the following results (6.1, 6.2 and 6.3 from \cite{kpw}).

\begin{theorem}
 Let $(X, \varphi)$ be a mixing Smale space and $P,Q$ be finite $\varphi$-invariants sets.
 \begin{enumerate}
  \item If $a$ is in $C_{c}(G^{s}(X, \varphi, Q))$ and $b$ is in $C_{c}(G^{u}(X, \varphi, P))$, then $ab$ and 
  $ba$ are both finite rank operators on $l^{2}(X^{h}(P,Q))$.
  \item If $a$ is in $\SX$ and $b$ is in $\UX$, then $ab$ and 
  $ba$ are both compact operators on $l^{2}(X^{h}(P,Q))$.
  \item Assuming that $P \cap Q = \emptyset$, if $a$ is in $\SX$ and $b$ is in $\UX$, then
  \[
   \lim_{n \rightarrow + \infty} \alpha^{-n}(a) b \, =  \lim_{n \rightarrow + \infty} b \alpha^{-n}(a)= 0.
  \]
\item If $a$ is in $\SX$ and $b$ is in $\UX$, then
  \[
   \lim_{n \rightarrow + \infty} \parallel \alpha^{n}(a) \alpha^{-n}(b) - \alpha^{-n}(b) \alpha^{n}(a) \parallel= 0.
  \]
 \end{enumerate}
\end{theorem}

Our main objective in this paper is to study the asymptotic behavior of $Tr(\alpha^{n}(a) \alpha^{-n}(b))$, for
$a$ in $C_{c}(G^{s}(X, \varphi, Q))$ and $b$ in $C_{c}(G^{u}(X, \varphi, P))$ and relate it to $\tau^{s}(a)$ and $\tau^{u}(b)$.

\section{Main Result}\label{results}

We consider the case that $(X, \varphi)$ is a \emph{mixing} Smale space and show how the above traces on $\SX$ and $\UX$ are related to asymptotics of the usual trace on ${\mathfrak B}(l^2(X^h(P,Q)))$.  That is, for an 
operator $A$ on $l^2(X^h(P,Q))$, we have 
$$
Tr(A) = \sum_{x \in X^h(P,Q)}{<A\delta_x, \delta_x>}
$$
This is defined for trace class operators and for positive operators, allowing $+\infty$
as a possible value.

\begin{theorem} \label{AsympTrace}
Let $(X, \varphi)$ be a mixing Smale space with topological entropy $\log(\lambda)$, and let $a \in \SX$, $b \in \UX$. If either
\begin{enumerate}
\item $a \in C_c(\GsX)$ and $b \in C_c(\GuX)$, or
\item $a$ and $b$ are both positive
\end{enumerate}
then
$$
\lim_{k \rightarrow +\infty}\lambda^{-2k} Tr\left(\alpha^k(a)\alpha^{-k}(b) \right) = \tau^s(a)\tau^u(b).
$$
In the second case, both sides of the equality may be $+\infty$.
\end{theorem}

To begin the proof, we must describe the topologies on $G^{s}(X, \varphi, Q)$ and $G^{u}(X, \varphi, P)$.
We take the following from  \cite{putnam96} or  \cite{whittaker11}.
To consider the former, suppose that $(x,y)$ is in $G^{s}(X, \varphi, Q)$.
It follows that for some positive integer $n$, $\varphi^{n}(y) \in X^{s}(\varphi^{n}(x, \epsilon_{X}/2)$.
Then the function
\[
 h^{s}(z) = \varphi^{-n}[  \varphi^{n}(z), \varphi^{n}(x)]
\]
sends $y$ to $x$. It is also defined on $X^{u}(y, \delta)$ for some $\delta > 0$ and is a 
homeomorphism to its image in $X^{u}(x, \epsilon_{X})$, which is open. Finally, we have 
that $h^{s}(z)$ is stably equivalent to $z$ for all $z$ in the domain. We define
\[
 X^s(x,y,\delta) = \{ (h^{s}(z), z) \mid z \in X^{u}(y, \delta), h^{s}(z) \in X^{u}(x, \delta) \}.
\]
We have
$(x, y) \in  X^s(x,y,\delta) \subset G^{s}(X, \varphi, Q)$. 
Such sets (allowing $x, y, \delta$ to vary)
form a neighbourhood base
for the topology of $G^{s}(X, \varphi, Q)$. The topology on 
$G^{u}(X, \varphi, P)$ is defined in an analogous way via 
functions denoted $h^{u}$ and sets $X^{u}(x, y, \delta)$, respectively.

The representation of $\SX$ is as follows, the representation of $\UX$ is completely analogous. 
Let $a$ be a function with compact support in $X^s(x,y,\delta)$, so that
$a$ is in $C_{c}(G^{s}(X, \varphi, Q)$. In fact, every element of $C_{c}(G^{s}(X, \varphi, Q)$
is a sum of such functions. For $w \in X^h(P,Q)$, we have 
$$
a\delta_w(x) = \sum_{(x,y) \in G^s(X,\varphi, Q)} a(x,y)\delta_w(y) = \left\{ \begin{array}{cc}
a(h^s(w),w)\delta_{h^s(w)} & \textrm{if} \ w \in X^u(y,\delta), \ h^s(w) \in X^u(x,\delta) \\
0 & \textrm{if} \ w \notin X^u(y,\delta)
\end{array} \right.
$$

For convenience, we now extend the class of functions that we consider on $\GsX$, $\GuX$.  As in Chapter II of \cite{renault} we consider the set of bounded Borel functions with compact support on these groupoids, denoted $B(\GsX)$ ($B(\GuX)$ respectively).  With convolution and involution defined as for the continuous functions, $B(\GsX)$ is a $\ast$-algebra.  Furthermore, every representation of $C_c(\GsX)$ extends to a representation of $B(\GsX)$.  This allows us to consider $Tr(a)$ where $a \in B(\GsX)$.  Moreover, since $B(\GsX)$ consists of integrable functions, we can make the following definition.

\begin{definition}
For $a \in B(G^{s}(X, \varphi,Q))$ and $b \in B(G^{u}(X, \varphi, P))$, define
\begin{eqnarray*}
% \tau^h(f) &=& \int_X f(x,x) d\mu, \\
I^s(a) &=&  \int_{X^u(Q)} a(x,x) d\mu^u, \ \textrm{and} \\
I^u(b) &=& \int_{X^s(P)} b(x,x) d\mu^s.
\end{eqnarray*}
\end{definition}

We start with the following lemma.
\begin{lemma}\label{0Trace}
Let $(X,\varphi)$ be a mixing Smale space. Let $a \in B(\GsX)$ be supported on a set of the 
form $X^s(x_a,y_a,\delta_a)$, and $b \in B(\GuX)$ supported on $X^u(x_b,y_b,\delta_b)$. We use $h^{s}_{a}$ 
and $h^{u}_{b}$ to denote the two functions as above.
Let $X_a = r(X(x_a,y_a,\delta_a))$, $X_b = s(X(x_b,y_b,\delta_b))$, where $r,s$ are the two canonical maps from 
$G^{s}(X, \varphi, Q)$ to $X^{u}(Q)$. 
If $h^s_a = id$ and $h^u_b = id$, then for each $k \in \mathbb{N}$
$$
Tr(\alpha^k(a)\alpha^{-k}(b)) = \sum_{w \in \varphi^k(X_a) \cap \varphi^{-k}(X_b)}a(\varphi^{-k}(w),\varphi^{-k}(w))b(\varphi^{k}(w),\varphi^{k}(w))
$$
otherwise
$$
\lim_{k \rightarrow +\infty}Tr(\alpha^k(a)\alpha^{-k}(b)) = 0
$$
\end{lemma}
\paragraph{Proof:}
If $x \neq h^s_a(y)$, $a(x,y) = 0$. to simplify notation, let $a(y) = a(h^s_a(y),y)$. Similarly, let $b(y) = b(h^u_b(y),y)$.
Now
\begin{eqnarray*}
Tr(\alpha^k(a)\alpha^{-k}(b)) &=& \sum_{w \in X^h(P)}<\alpha^k(a)\alpha^{-k}(b)\delta_w,\delta_w> \\
&=& \sum_{w \in X^h(P)} a(\varphi^{-2k}h^u_b\varphi^k(w)) b(\varphi^k(w))<\delta_{\varphi^kh^s_a\varphi^{-2k}h^u_b\varphi^k(w)}, \delta_w>.
\end{eqnarray*}
Suppose $h^s_a = h^u_b = id$, then we have
$$
Tr(\alpha^k(a)\alpha^{-k}(b))  =
\sum_{w \in X^h(P)}a(\varphi^{-k}(w))b(\varphi^{k}(w))<\delta_w,\delta_w>.
$$
Moreover, $a(\varphi^{-k}(w)) \neq 0$ only if $w \in \varphi^k(X_a)$, and $b(\varphi^{k}(w)) \neq 0$ only if $w \in \varphi^{-k}(X_b)$, so
$$
Tr(\alpha^k(a)\alpha^{-k}(b)) = \sum_{w \in \varphi^k(X_a) \cap \varphi^{-k}(X_b)}a(\varphi^{-k}(w))b(\varphi^{k}(w)).
$$

Now suppose $h^s_a \neq id$.
$$
Tr(\alpha^k(a)\alpha^{-k}(b)) = 
\sum_{w \in E_k} a(\varphi^{-2k}h^u_b\varphi^k(w)) 
b(\varphi^k(w))<\delta_{\varphi^kh^s_a\varphi^{-2k}h^u_b\varphi^k(w)}, \delta_w>,
$$
where 
$$
E_k = \{w \in X^h(P) | \ \varphi^k(w) \in X_b, \ \varphi^{-2k}h^u_b\varphi^k(w) \in X_a,  \ \varphi^kh^s_a\varphi^{-2k}h^u_b\varphi^k(w) = w   \}.
$$

We show that for large enough $k$, the set $E_k$ is empty.  We can find $\delta>0$ such that 
$$
d(z,h^s_a(z))>\delta
$$
for all $z \in X_a$. Now, we can find $K$ sufficiently large so that for all $k >K$ we have
$$
d(\varphi^{-2k}(y),\varphi^{-2k}h^u_b(y)) < \delta
$$
for all $y \in X_b$. So for $w \in E_k$, $ \varphi^k(w) \in X_b$ and thus
$$
d(\varphi^{-k}(w),\varphi^{-2k}h^u_b(\varphi^k(w))) < \delta
$$
but $\varphi^{-2k}h^u_b(\varphi^k(w)) \in X_a$ so
$$
d(\varphi^{-2k}h^u_b(\varphi^k(w)), h^s_a\varphi^{-2k}h^u_b(\varphi^k(w))) > \delta
$$
so $\varphi^{-k}(w) \neq h^s_a\varphi^{-2k}h^u_b(\varphi^k(w))$, contradicting $w \in E_k$. Hence for large $k$, $E_k$ is empty and 
$$
Tr(\alpha^k(a)\alpha^{-k}(b)) = 0.
$$
A similar argument gives the result in the case $h^u_b \neq id$.
\qed

\begin{lemma}\label{stepresult}
Let $(X,\varphi)$ be a mixing Smale space. Let $a \in B(\GsX)$  and $b \in B(\GuX)$ be step functions. Then 
$$
\lim_{k\rightarrow \infty}\lambda^{-2k}Tr(\alpha^k(a)\alpha^{-k}(b)) = I^s(a)I^u(b).
$$
\end{lemma}

\paragraph{Proof:}
By Lemma \ref{0Trace} we can assume, without loss of generality that 
$a(x,y) = 0$ whenever $x \neq y$, and $b(x,y) = 0$ if $x \neq y$.  
We therefore consider $a(x,x)$ to be supported on $A \subset X^s(x_a)$, and similarly with $b$.
 
Let $\{A_{j}\}_{j=1}^{\infty}$ be precompact sets on which $a$ is constant, and $a_{j}$ the value on $A_{j}$.  Similarly $\{B_{j}\}_{j=1}^{\infty}$ and $b_{j}$.  Then for each $k$ we have

\begin{eqnarray*}
I^s(a) &=& \sum_{j=1}^{\infty}\mu^u(A_{j})a_{j} \\
I^u(b) &=& \sum_{j=1}^{\infty}\mu^u(B_{j})b_{j} \\
Tr(\alpha^k(a)\alpha^{-k}(b)) &=& \sum_j\sum_l \#h^k_{j,l}a_{j}b_{l}
\end{eqnarray*}
Where $h^k_{j,l} = \varphi^k(A_{j}) \cap \varphi^{-k}(B_{l})$.  From Theorem 2.4 of \cite{killoughputnam11} we know that
$$
\lim_{k \rightarrow \infty} \lambda^{-2k}h^k_{j,l} = \mu^u(A_{j})\mu^s(B_{l}).
$$

So
\begin{eqnarray*}
\lim_{k\rightarrow \infty}\lambda^{-2k}Tr(\alpha^k(a)\alpha^{-k}(b)) &=& \sum_j\sum_l \lim_{k \rightarrow \infty} \lambda^{-2k}\#h^k_{j,l}a_{j}b_{l}\\
&=& \sum_j\sum_l\mu^u(A_{j})\mu^s(B_{l})a_{j}b_{l} \\
&=& (\sum_j\mu^u(A_{j})a_{j})(\sum_l\mu^s(B_{l})b_{l}) \\
&=& I^s(a)I^u(b).
\end{eqnarray*}
\qed

\paragraph{Proof of Theorem \ref{AsympTrace}:}
We start with the special case that $a \in C_c(\GsX)$ and $b \in C_c(\GuX)$ are both positive.  Lemma \ref{0Trace} allows us to assume, without loss of generality, that $a(x,y) = 0$ whenever $x \neq y$, and $b(x,y) = 0$ if $x \neq y$. 
%We can find increasing sequences of step functions $\{a_i\}_1^{\infty}$, $\{b_i\}_1^{\infty}$ such that $a_i \rightarrow a$ and $b_i \rightarrow b$. Then by monotone convergence $\tau^s(a_i) \rightarrow \tau^s(a)$, $\tau^u(b_i) \rightarrow \tau^u(b)$. Moreover, we also have, for each $k$
%$$
%lim_{i \rightarrow \infty} Tr(\alpha^k(a_i)\alpha^{-k}(b_i)) = Tr(\alpha^k(a)\alpha^{-k}(b))
%$$

%Let $\{A_{i,j}\}_{j=1}^{\infty}$ be sets on which $a_i$ is constant, and $a_{i,j}$ the value on $A_{i,j}$.  Similarly $\{B_{i,j}\}_{j=1}^{\infty}$ and $b_{i,j}$.  Then for each $k$ we have

%\begin{eqnarray*}
%\tau^s(a_i) &=& \sum_{j=1}^{\infty}\mu^u(A_{i,j})a_{i,j} \\
%\tau^u(b_i) &=& \sum_{j=1}^{\infty}\mu^u(B_{i,j})b_{i,j} \\
%Tr(\alpha^k(a_i)\alpha^{-k}(b_i)) &=& \sum_j\sum_l \#h^k_{j,l}a_{i,j}b_{i,l}
%\end{eqnarray*}
%Where $h^k_{j,l} = \varphi^k(A_{i,j}) \cap \varphi^{-k}(B_{i,l})$.  From Theorem 2.3 in \cite{killoughputnam11} we know that
%$$
%\lim_{k \rightarrow \infty} \lambda^{-2k}h^k_{j,l} = \mu^u(A_{i,j})\mu^s(B_{i,l}).
%$$

%So
%\begin{eqnarray*}
%\lim_{k\rightarrow \infty}\lambda^{-2k}Tr(\alpha^k(a_i)\alpha^{-k}(b_i)) &=& \sum_j\sum_l \lim_{k \rightarrow \infty} \lambda^{-2k}\#h^k_{j,l}a_{i,j}b_{i,l}\\
%&=& \sum_j\sum_l\mu^u(A_{i,j})\mu^s(B_{i,l})a_{i,j}b_{i,l} \\
%&=& (\sum_j\mu^u(A_{i,j})a_{i,j})(\sum_l\mu^s(B_{i,l})b_{i,l}) \\
%&=& \tau^s(a_i)\tau^u(b_i).
%\end{eqnarray*}

%Taking limits as $i \to \infty$ we have
%$$
%\lim_{k \rightarrow +\infty}\lambda^{-2k} Tr\left(\alpha^k(a)\alpha^{-k}(b) \right) = \tau^S(a)\tau^U(b).
%$$
We find bounded Borel step functions with compact support $a_l, a_u$ and $b_l, b_u$ such that 
\begin{eqnarray*}
a_l \leq a \leq a_u, &  & b_l \leq b \leq b_u, \\
I^s(a_l) + \epsilon > \tau^s(a) > I^s(a_u)- \epsilon, & \textrm{and} & I^u(b_l) + \epsilon > \tau^u(b) > I^u(b_u)- \epsilon.
\end{eqnarray*}

For each $k$
$$
\lambda^{-2k}Tr\left(\alpha^{k}(a_l)\alpha^{-k}(b_l)\right) \leq \lambda^{-2k}Tr\left(\alpha^{k}(a)\alpha^{-k}(b)\right) \leq \lambda^{-2k}Tr\left(\alpha^{k}(a_u)\alpha^{-k}(b_u)\right).
$$
So, by Lemma \ref{stepresult}, in the limit as $k \to \infty$ we have
$$
I^s(a_l)I^u(b_l) \leq \lim _{k \to \infty} \lambda^{-2k}Tr\left(\alpha^{k}(a)\alpha^{-k}(b)\right) \leq I^s(a_u)I^u(b_u).
$$

Therefore, for all $\epsilon >0$, we have
$$
\left( \tau^s(a) - \epsilon \right)\left( \tau^u(b) - \epsilon \right) <  \lim _{k \to \infty} \lambda^{-2k}Tr\left(\alpha^{k}(a)\alpha^{-k}(b)\right) <\left( \tau^s(a) + \epsilon \right)\left( \tau^u(b) + \epsilon \right).
$$
Hence, we conclude
$$
\lim _{k \to \infty} \lambda^{-2k}Tr\left(\alpha^{k}(a)\alpha^{-k}(b)\right) = \tau^s(a)\tau^u(b).
$$

To prove case (1) of the theorem (the $C_c$ case) we simply 
apply the above argument separately to the positive/negative parts of the real/imaginary components of $a$, $b$.

Now suppose $a$, $b$ are positive.
If $\tau^s(a)$ and $\tau^u(b)$ are finite, then as above we 
can approximate above and below with (integrable) step functions, and by
Lemma \ref{stepresult}, the theorem holds. Now suppose one of the traces
is not finite.  Without loss of generality,  assume $\tau^s(a) = +\infty$. We can then find $\{a_i\}$ such that $a_i \leq a$, $a_i$ increasing, $a_i \in C_c(\GsX)$ positive, and $\tau^s(a_i) \geq i$. Similarly choose $\{b_i\}$ (if $b \in C_c(\GuX)$, choose $b_i = b$ for all $i$). Then we have

\begin{eqnarray*}
\lim_{k\rightarrow \infty}\lambda^{-2k}Tr(\alpha^k(a)\alpha^{-k}(b)) & \geq & \lim_{k\rightarrow \infty}\lambda^{-2k}Tr(\alpha^k(a_i)\alpha^{-k}(b_i)) \\
 &=& \tau^s(a_i)\tau^u(b_i) \\
 & \geq & i\tau^u(b_1).
\end{eqnarray*}

As this holds for any $i \in \mathbb{Z}$, we have
$$
\lim_{k\rightarrow \infty}\lambda^{-2k}Tr(\alpha^k(a)\alpha^{-k}(b)) = \infty = \tau^s(a)\tau^u(b).
$$
\qed

%	\TOCadd{Bibliography}
	\bibliographystyle{plain}
	\bibliography{references}

\end{document}